\documentclass[conference,final,a4paper]{IEEEtran}

\makeatletter
\def\ps@IEEEtitlepagestyle{%
  \def\@oddfoot{\mycopyrightnotice}%
  \def\@evenfoot{}%
}
\def\mycopyrightnotice{%
  {\footnotesize \copyright~2022 IEEE\hfill}
  \gdef\mycopyrightnotice{}
}

\usepackage{blindtext}
\usepackage{eso-pic}
\usepackage{cite}
\usepackage{amsmath,amssymb,amsfonts}
\usepackage{algorithmic}
\usepackage{graphicx}
\usepackage{url}
\usepackage{orcidlink}
\usepackage{siunitx}
\usepackage{textcomp}
\usepackage{booktabs}
\usepackage{xcolor}
\usepackage{hyperref}
\def\BibTeX{{\rm B\kern-.05em{\sc i\kern-.025em b}\kern-.08em
    T\kern-.1667em\lower.7ex\hbox{E}\kern-.125emX}}
    
\usepackage{eso-pic}

\begin{document}
\title{\vspace*{1cm} 
Efficient Simulation of Complex Capillary Effects in Advanced Manufacturing Processes using the Finite Volume Method
}

\author{\IEEEauthorblockN{Patrick Zimbrod\orcidlink{0000-0003-3108-3171}}
\IEEEauthorblockA{\textit{Digital Manufacturing} \\
\textit{University of Augsburg}\\
Augsburg, Germany\\
Patrick.Zimbrod@uni-a.de}
\and
\IEEEauthorblockN{Magdalena Schreter\orcidlink{0000-0003-3888-4086}}
\IEEEauthorblockA{\textit{Computational Mechanics} \\
\textit{Technical University Munich}\\
Munich, Germany\\
Magdalena.Schreter@tum.de}
\and
\IEEEauthorblockN{Johannes Schilp}
\IEEEauthorblockA{\textit{Digital Manufacturing} \\
\textit{University of Augsburg}\\
Augsburg, Germany\\
Johannes.Schilp@uni-a.de}
}

\maketitle
\begin{abstract}
The accurate representation of surface tension driven flows in multiphase systems is considered a challenging problem to resolve numerically. Although there have been extensive works in the past that have presented approaches to resolve these so called Marangoni flows at the phase boundaries, the question of how to efficiently resolve the interface in a universal and conservative manner remains largely open in comparison. Such problems are of high practical relevance in many manufacturing processes, especially in the microfluidic regime where capillary effects dominate the local force equilibria. In this work, we present a freely available numerical solver based on the Finite Volume Method that is able to resolve arbitrarily complex, incompressible multiphase systems with the mentioned physics at phase boundaries. An efficient solution with respect to the number of degrees of freedom can be obtained by either using high order WENO stencils or by employing adaptive cell refinement. We demonstrate the capabilities of the solver by investigating a model benchmark case as well as a single track laser melting process that is highly relevant within laser additive manufacturing.
\end{abstract}


\section{Introduction}

Accurate modelling of complex physics has received considerable attention within the last two decades.
This can mostly be attributed to the fact that an increased understanding of the phenomena relevant for manufacturing directly leads to more finely tuned or even new processes that tend to increase output quality.

Two such examples that shall serve as a guideline for this work are microfluidic applications and additive manufacturing of metals.
Both groups of processes are characterized by complex flow phenomena that involve several physical phenomena on a microscopic scale \cite{sackmannPresentFutureRole2014,debroyAdditiveManufacturingMetallic2018}. The fact that those effects tend to have fast dynamics further complicates the analysis \cite{khairallahControllingInterdependentMesonanosecond2020}. This leads to an increased amount of empirical effort necessary to capture and quantify the flow patterns involved.

Hence, simulation has become an important alternative to extensive experimental research. Within computational simulation frameworks, the involved physics can be precisely monitored, even down to a nanoseconds scale, which would oftentimes be probhibitively expensive to do experimentally. At the same time, a considerable drawback of this approach is the additional amount of work that has to be put into the mathematical modelling of the involved physics. Especially within fluid dynamics, gaining a stable as well as accurate solution can oftentimes be challenging, as the resulting set of partial differential equations can have a small stability region in the temporal domain \cite{laxStabilityDifferenceSchemes2013}.

Therefore in this work, we propose and demonstrate a performant, free and open source software framework based on the Finite Volume method that is able to capture complex, non-isothermal surface flows. We discuss the necessary mathematical modelling in order to extend the existing solver ecosystem towards incorporating non-isothermal surface tension driven flows. We showcase the capabilities of the developed application using two example cases with different geometric and physical complexity.

\section{Related works}

Modelling field-dependent capillary convection, also known as Marangoni convection, in itself is not a new field of research.
Various works exist that address the need to develop robust numerical methods in order to resolve the complex flow patterns that arise. Those are mostly driven by application, i.e. they are tailored towards simulating a particular set of physics.

Early works include the investigation of marangoni effects in microgravity \cite{chunExperimentsTransitionSteady1979,chunMicrogravitySimulationMarangoni1978}, buoyancy-driven flows \cite{villersCoupledBuoyancyMarangoni1992} and welding \cite{limmaneevichitrExperimentsSimulateEffect2000,limmaneevichitrVisualizationMarangoniConvection2000}.

After that, Marangoni convection in the microfluidic regime has received increasing attention within the last ten years \cite{karbalaeiThermocapillarityMicrofluidicsReview2016}. Numerous scenarios have been investigated since, including engineering of traps and pumps, stability in microchannels and particle accumulation in microfluidic flows \cite{panInstabilityMarangoniToroidal2011,basuVirtualMicrofluidicTraps2008,orlishausenParticleAccumulationDepletion2017}. More recent works have also put a strong focus towards increasing the performance of such models substantially \cite{kronbichler2018fast}.

Another field where the effects of capillary convection have become a dominant part of the involved physics is within metal additive manufacturing. Here, the focus lays on resolving temperature-dependent surface tension flows, as melt pools exhibit large spatial and temporal temperature gradients \cite{debroyAdditiveManufacturingMetallic2018}. Various methods have been employed to resolve the thermo-fluid dynamics on the powder scale including Arbitrary Lagrangian Eulerian \cite{khairallahMesoscopicSimulationModel2014,martinDynamicsPoreFormation2019}, Finite Difference \cite{gusarovModelRadiationHeat2009}, Finite Element \cite{zhangMultiscaleMultiphysicsModeling2018,caoNovelHighefficientFinite2021,meier2021physics}, Lattice Boltzmann \cite{krzyzanowskiMultiphysicsSimulationApproach2021}, Smoothed Particle Hydrodynamics \cite{wimmerExperimentalNumericalInvestigations2021} and the Finite Volume Method which shall be used for this work as well \cite{leeMesoscopicSimulationHeat2015,gurtlerSimulationLaserBeam2013,ottoMultiphysicalSimulationLaser2012}.

However, the mentioned works exhibit some shortcomings when compared to the general current state of numerical modelling. Especially in applications where highly dynamic physics are strongly localized in the computational domain, adaptively refining the domain is a very desirable property that is not widely employed yet. Furthermore, the use of high order spatial approximations remains largely restricted to few numerical schemes of Galerkin type. However, achieving high order accuracy can help alleviate the need for overly fine computational meshes and hence is also desirable. Lastly, few of the presented solution along with their source codes are freely available. To the authors' best knowledge, none of the current works on the presented topic incorporate all three properties. Thus, we present a framework that allows for grid adaptivity, high order schemes and is openly available (see section \ref{sec:data} for references).

\section{Theory}

To describe the underlying mathematical concepts, we will use the index notation of tensors in a cartesian coordinate system alongside with the Einstein summation convention where the sum over matching indices is implied.

From here on, we assume a system of partial differential equations that involves the transient evolution of multiple fluid phases in a joint domain. This means we are tasked to solve at least the Navier-Stokes equations in order to capture velocity, pressure and temperature fields.

In addition, we need to properly discretize the different phases in the model. One way to achieve this in a Finite Volume framework is through the Volume-of-Fluid method \cite{hirtVolumeFluidVOF1981}:

\begin{equation}
  \frac{\partial \alpha_{j}}{\partial t} + \frac{\partial u_{i} \alpha_{j}}{\partial x_{i}} = 0
\end{equation}

Where $\alpha_{j}$ represents the $j$-th phase involved. From this, we obtain a set of hyperbolic conservation laws that track the cell-wise volumetric fraction of each phase over the domain.

\begin{figure}[!tbp]
  \centering
  \includegraphics[width=\linewidth]{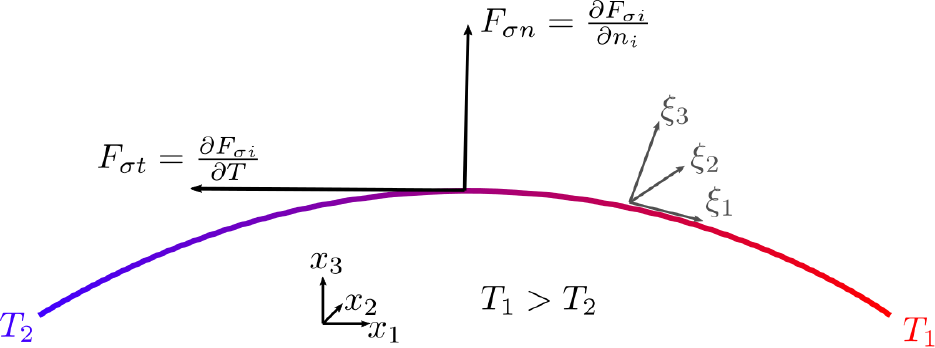}
  \caption{Forces acting on an interface between two immiscible phases subject to a temperature gradient $\Delta T = T_{2} - T_{1}$. The overall surface tension Force $F_{\sigma}$ can be divided into a normal component $F_{\sigma n}$ and a tangential component $F_{\sigma t}$ in the coordinate system $\xi_{i}$ proprietary to the interface. The tangential component is also called Marangoni force and depends on the gradient of external fields.}
  \label{fig:surface-tension}
\end{figure}

After that, we can consider the effects happening at the interface of such phases, i.e. within this work we investigate surface tension. With no other physics involved, this leads to a contractional movement normal to the interface. If however varying fields, such as temperature are present that influence the surface tension coefficient, additional forces interact with the interface. This is shown in Fig. \ref{fig:surface-tension}.

The theory of the proposed implementation of those surface tension effects relies on the Continuum Surface Stress Method initially presented by Lafaurie et al. \cite{lafaurieModellingMergingFragmentation1994}, extending the works on modelling surface tension with the Volume of Fluid Method as described by Brackbill et al. \cite{brackbillContinuumMethodModeling1992}. Using this modelling technique, it is possible to incorporate surface tension as a generalized body force into the governing equations that is well defined everywhere in the domain. This additional term is then called the capillary stress which takes the form of a stress tensor:

\begin{equation}
    T_{ij} = - \sigma \delta_s (\delta_{ij} - n_i n_j)
\end{equation}

Where $\sigma$ is the generally non constant coefficient of surface tension. $\delta_s$ is the interface delta function that serves as an indicator function of where the phase boundaries are located. $\delta_{ij}$ is the Kronecker Delta function and $n_i$ is the interface unit normal vector. The concrete implementation of these abstract quantities into the Volume of Fluid framework is given later on. The local coordinate system of the interface needs to be an orthogonal system in order to separate the purely geometry based capillary force and the tangential Marangoni forces (c.f. Fig. \ref{fig:surface-tension}).
The resulting force accounting for all surface tension effects can now be expressed as the negative divergence of the capillary stress tensor, yielding:

\begin{multline} \label{eq:stress-tensor}
    \frac{\partial T_{ij}}{\partial x_j} = \frac{\partial \sigma}{\partial x_j} [\delta_s(\delta_{ij} - n_i n_j)] + \frac{\partial \delta_s}{\partial x_j} [\sigma (\delta_{ij} - n_i n_j)] \\ + \sigma \delta_s \frac{\partial}{\partial x_j} (\delta_{ij} - n_i n_j)
\end{multline}

One can show that Eq. \ref{eq:stress-tensor} can be re-arranged into a much shorter and more useful form. The derivation is given in more detail by Lafaurie et al. \cite{lafaurieModellingMergingFragmentation1994}. By performing those rearrangements, we can recover a form that separates the normal from the tangential components of the divergence vector:

\begin{equation}
\label{eq:divergence-capillarystress}
    \frac{\partial T_{ij}}{\partial x_j} = \frac{\partial \sigma}{\partial x_i} [\delta_s(\delta_{ij} - n_i n_j)] - \sigma \kappa n_i \delta_s 
\end{equation}

Where we introduced the interface curvature $\kappa$. The second term of the right hand side corresponds to the normal capillary force (Fig. \ref{fig:surface-tension}, $F_{\sigma n}$) directed in the normal direction of the interface, effecting a contractional movement of the interface. The first term resembles the marangoni-type forces (Fig. \ref{fig:surface-tension}, $F_{\sigma t}$) present. By taking the derivative of the surface tension coefficient, this term does not vanish if and only if spatial gradients exist at the interface. This is normally the case when there are multiple species involved or temperature gradients present \cite{j.straubThermokapillareGrenzflachenkonvektionGasblasen1990}. Note that we otherwise simply evaluate the ordinary form of surface tension without any additional physics present. However, as the surface tension coefficient is otherwise not a direct local variable, we must further differentiate the term:

\begin{equation}
    \frac{\partial \sigma(c,T)}{\partial x_j} = \frac{\partial \sigma(c,T)}{\partial T} \frac{\partial T}{\partial x_j} + \frac{\partial \sigma(c,T)}{\partial c} \frac{\partial c}{\partial x_j}
    \label{eq:sigmadiff}
\end{equation}

Additionally, the delta function $\delta_s$ still needs to be discretized in a suitable manner in order to capture the physics within the Finite Volume framework. This can be done using the Volume of Fluid Method by taking the gradient of the phase volume fraction $\alpha$ \cite{gueyffierVolumeofFluidInterfaceTracking1999,hirtVolumeFluidVOF1981}:

\begin{equation}
    \delta_s = \left\lvert \frac{\partial \alpha}{\partial x_i} \right\rvert
\end{equation}

The interface unit normal vector $n_i$ can be computed in a similar manner by using the previously computed interface function:

\begin{equation}
    n_i = \frac{1}{\delta_s} \frac{\partial \alpha}{\partial x_i}
\end{equation}


\section{Implementation}

The presented theoretical foundation is modified and implemented based on the open source library \textit{OpenFOAM}, which is a popular software library that implements the Finite Volume Method \cite{jasakOpenFOAMLibraryComplex2007}. The modifications rely on the works of Gueyffier et al. in order to account for the additional Marangoni stresses \cite{gueyffierVolumeofFluidInterfaceTracking1999}.

It is clear from equation \ref{eq:sigmadiff} that we would normally have to supply the derivatives of the surface tension coefficient with respect to the respective fields and then compute the gradients of the fields each. This can quickly become tedious to implement as well to code when multiple sources of Marangoni convection are considered. Therefore, we instead first compute the surface tension coefficient locally and treat it as a differentiable field itself.
A very useful property of this approach is that not only thermally driven surface tension effects can be incorporated, but also gradients arising from any kind of inhomogeneity. This means that, among others, solute-driven effects can also be modelled. This scenario is common in manufacturing processes involving mixing of multiple solvents. Within the context of additive manufacturing, in-situ alloying is a scenario where such effects are not negligible \cite{katz-demyanetzInsituAlloyingNovel2020}.

In order to approximate the differential terms in the governing equations efficiently, we use high order weighted essentially non-oscillatory (WENO) schemes that have been developed separately for OpenFOAM by Gärtner et al. \cite{gartnerEfficientWENOLibrary2020,martinImplementationValidationSemiImplicit2018}. In this implementation, the numeric stencil coefficients to specified order are pre-computed and locally stored to alleviate the harsh overhead in memory requirements. The theory behind the derivation of high order WENO stencils is given in Gärtner et al. \cite{gartnerEfficientWENOLibrary2020}.

\section{Numerical experiments}

First, we replicate and examine the numerical setup proposed by Ma and Bothe that models the flow around a stationary bubble subject to a temperature gradient \cite{maDirectNumericalSimulation2011}. We aim to accurately predict the velocity field that arises due to marangoni convection at the interface.

We then proceed to investigate the mesoscopic nanosecond dynamics during selective laser melting of steel, a common problem within additive manufacturing. Here, we show that our implementation is able to handle geometrically and physically complex simulations in 3D.

\subsection{Marangoni Flow around a droplet}\label{sec:droplet}

We now investigate the aforementioned simple, two-dimensional benchmark case that illustrates the ability to capture Marangoni effects for the present work. An overview of the material parameters is given in table \ref{tab:droplet-parameters}.

\begin{table}[!t]
\renewcommand{\arraystretch}{1.3}
\caption{Material properties of the oscillating droplet case \cite{maDirectNumericalSimulation2011}. The indices indicate either a property of the surrounding medium (1) or of the droplet (2).}
\label{tab:droplet-parameters}
\centering
\begin{tabular}{ll}
\toprule
Quantity & Value [Unit]\\
\midrule
Density $\rho_{1}$ & \SI{250}{\kilo\gram \per \metre\cubed} \\
Density $\rho_{2}$ &  $2 \rho_{1}$ \\
Heat capacity $c_{p,1}$ & \SI{5e-5}{\joule \per \kilo\gram \per \kelvin} \\
Heat capacity $c_{p,2}$ & $2 c_{p,1}$ \\
Viscosity $\mu_{1}$ & \SI{0.012}{\kilo\gram \per \metre \per \second} \\
Viscosity $\mu_{2}$ & $2 \mu_{1}$ \\ 
Thermal conductivity $\lambda_{1}$ & \SI{1.2e-6}{\watt \per \metre \per \kelvin} \\
Thermal conductivity $\lambda_{2}$ & $2 \lambda_{1}$ \\
Surface tension coefficient $\sigma$ & \SI{0.1}{\newton \per \metre} \\
Marangoni coefficient $\sigma_{T}$ & \SI{0.02}{\newton \per \metre \per \kelvin} \\

\bottomrule
\end{tabular}
\end{table}

We initialize the droplet in the center of a rectangular domain of length $4a \times 4a$ where $a = \SI{1.44e-3}{\metre}$ with a radius of $a$. The droplet is subject to a temperature gradient of \SI{200}{\kelvin \per \metre} where the bottom temperature is fixed to \SI{290}{\kelvin}.

We discretize the domain using a coarse grid of 100 x 100 quadrilaterals. The spatial discretization of the divergence is done using a mixture of first order bounded Gaussian schemes and fourth order WENO schemes.
We impose no-slip boundary conditions at the walls as well as Neumann type boundaries for the temperature at the left and right walls. The temperature gradient is enforced during the entire simulation using Dirichlet boundaries at the bottom and top. We simulate the evolution of the involved fields up to $t = \SI{0.12}{\second}$.

The remainder of the involved numerics and parameters, as well as all boundary conditions can be accessed via the online repository given in section \ref{sec:data}.

\begin{figure}[!tbp]
  \centering
  \includegraphics[width=\linewidth]{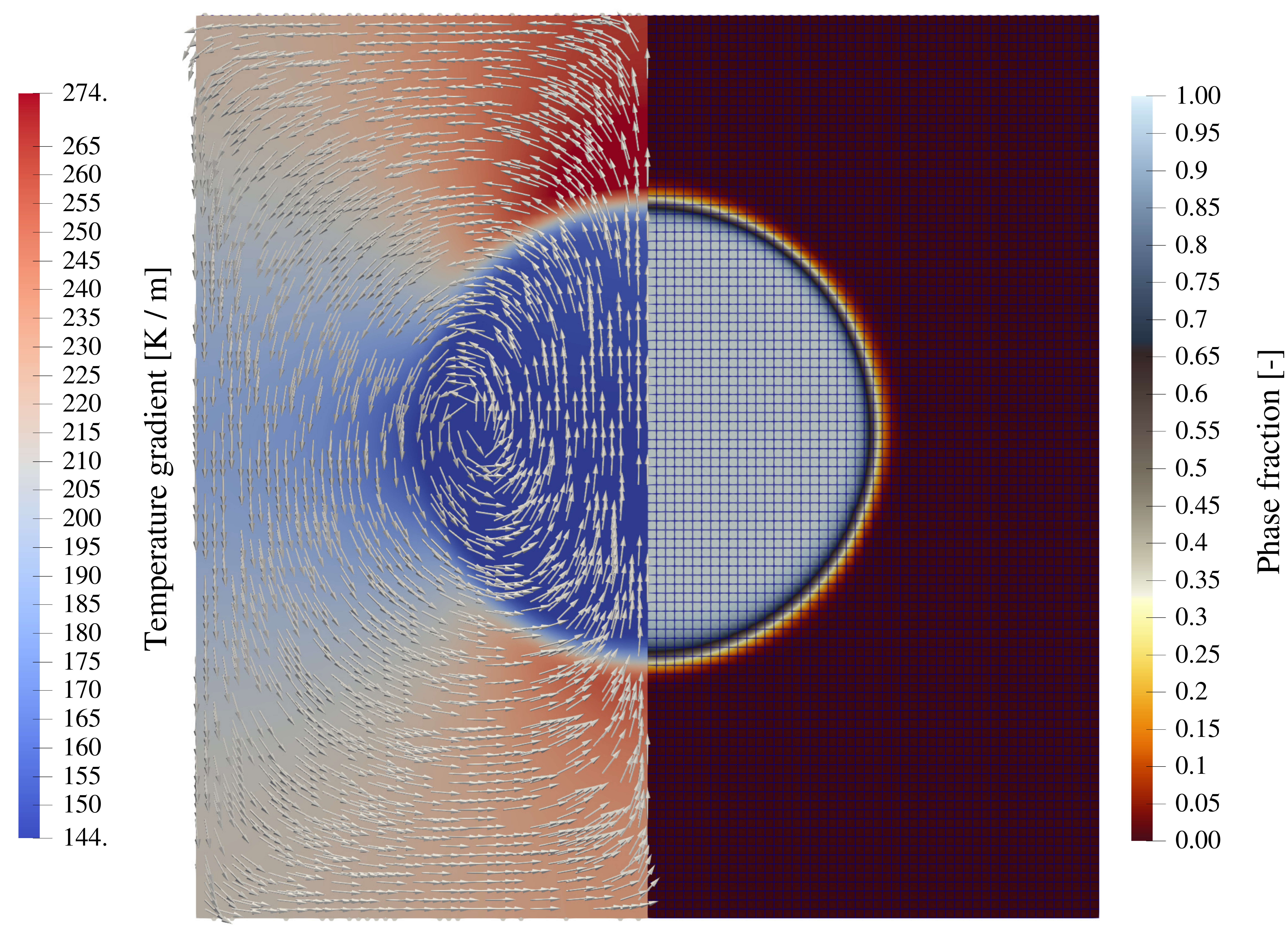}
  \caption{Flow field around the droplet at final time $t = \SI{0.12}{\second}$. The resulting velocity vector field (arrows) and temperature gradient field are shown on the left slice of the domain. The right side shows the phase fraction of the droplet as well as the spatial discretisation.}
  \label{fig:2d-benchmark}
\end{figure}

The resulting fields at the final time step $t = \SI{0.12}{\second}$ is shown in Fig. \ref{fig:2d-benchmark}. We can clearly observe the field of vortices that form around the bubble, initiated from varying pressures along the interface. This is due to the Marangoni effect, resulting in an uneven pressure distribution in the domain. Following equation \ref{eq:divergence-capillarystress}, the first term involving the gradient of the surface tension coefficient - i.e. the marangoni coefficient times the gradient of the temperature field - does not vanish here and hence induces a non homogeneous pressure contribution. The resulting two large vortices drive the droplet to an upwards motion in direction of increasing temperature gradient. The observed fields reproduce the findings of Ma and Bothe closely \cite{maDirectNumericalSimulation2011}.

\begin{figure}[!tbp]
  \centering
  \includegraphics[width=\linewidth]{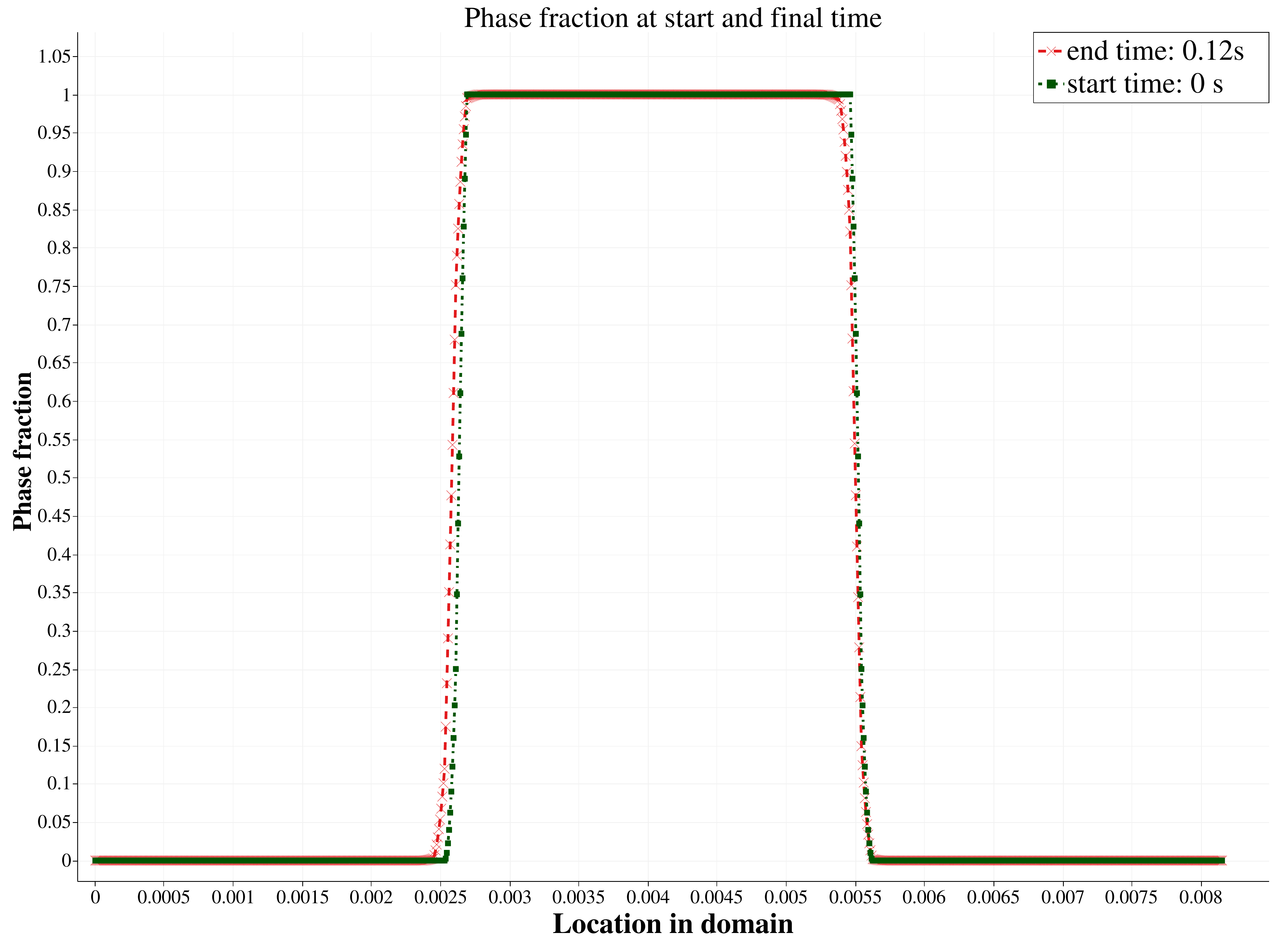}
  \caption{Phase fraction profiles of the droplet at beginning and end of the investigated temporal domain. This cross-section represents a diagonal slice from the left top to the right bottom of the grid shown in Fig. \ref{fig:2d-benchmark}.}
  \label{fig:alpha-comparison}
\end{figure}

Despite the comparatively low mesh resolution, the phase boundary is conserved reasonably well after solving for 6000 time steps. A comparison of the sharp interface between the two phases at the beginning and end of the simulation is given in Fig. \ref{fig:alpha-comparison}. During the solution process, only minimal numerical diffusion of the interface occurs, indicated by the slightly tapered off edges of the rectangular profile at end time. In addition, we don't observe any oscillatory behaviour in the fields which is essential for obtaining a stable solution.

\subsection{Powderbed-scale physics during selective laser melting}\label{sec:lpbf}

We now turn to a more complex and relevant simulation setup that appears regularly in additive manufacturing. We aim to investigate the melting, solidification and evaporation behavior of molten metal during the Laser Powder Bed Fusion Process (PBF-LB/M).

This case shall serve as a showcase of a physically complex problem that involves phase changes, simultaneous handling of solid, liquid and gaseous phases (and hence strong pressure gradients) as well as spatially and temporally localized heat transfer.

\begin{figure}[!tbp]
  \centering
  \includegraphics[width=\linewidth]{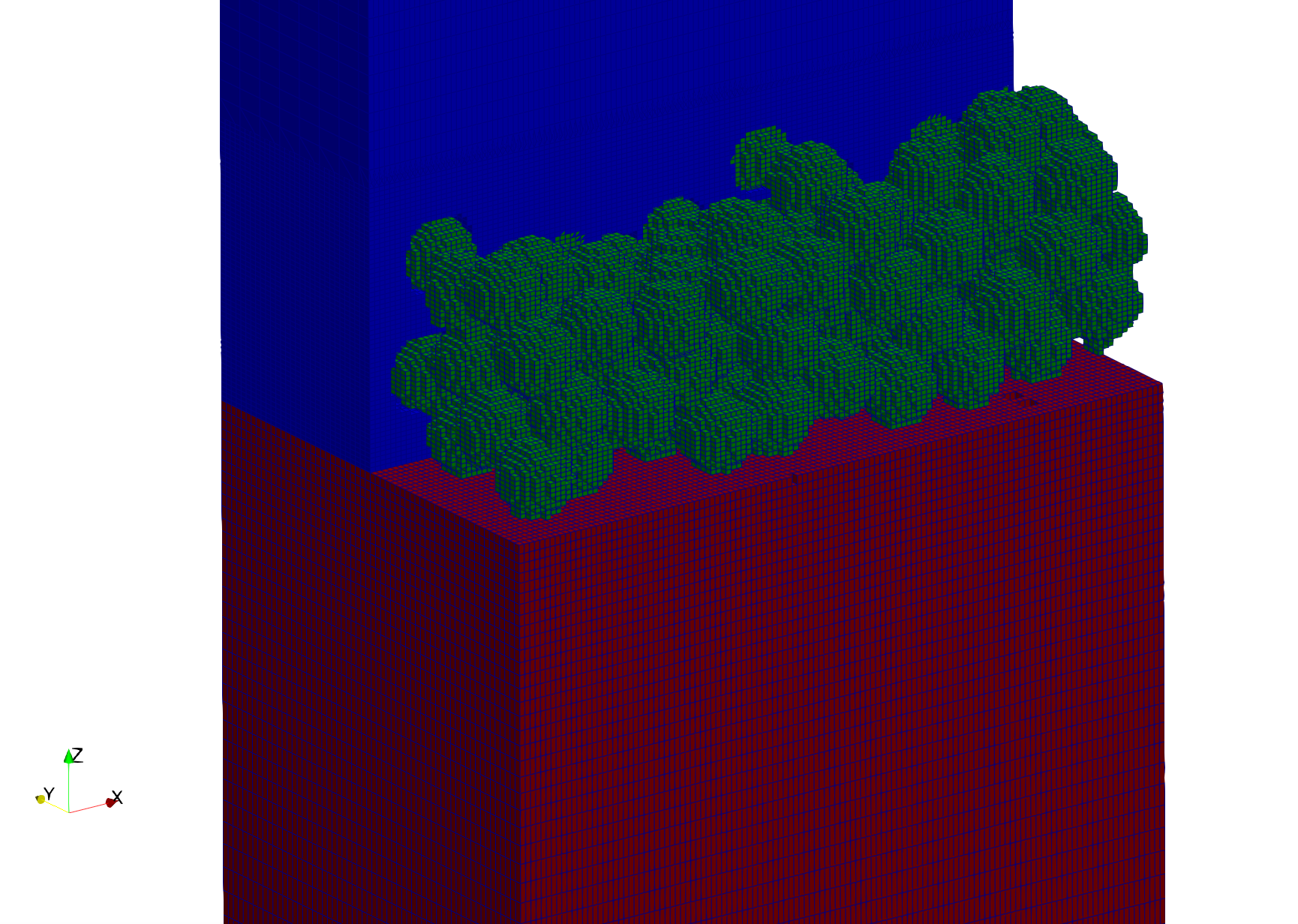}
  \caption{Composition of the computational domain at the initial time step. Solidified build platform (red), discretised powder bed in solid state (green) and gaseous argon atmosphere (blue). As reference, the individual cell boundaries are marked with blue lines. The initial grid is finest where the temperature gradients are expected to be highest.}
  \label{fig:lpbf-initial}
\end{figure}

The computational domain with the initial conditions for the phase fractions is shown in Fig. \ref{fig:lpbf-initial}. We create the packing of the metallic powder particles using a rain drop model implemented in the open source Discrete Element Method software Yade \cite{vaclavsmilauerYadeDocumentation2021}. Overall, the domain spans a cuboid geometry of size $\SI{0.5}{\milli\metre} \times \SI{0.3}{\milli\metre} \times \SI{0.7}{\milli\metre}$ and contains roughly 800.000 hexahedral cells arranged in a cartesian grid. The inhomogeneous temperature field is generated by a laser heat source of diameter \SI{50}{\micro\metre} that inputs an amount of heat equivalent to \SI{67}{\joule \per \milli\metre\cubed}.

The physical parameters of the PBF-LB/M model are given in table \ref{tab:lpbf-parameters}. The numerical parameters are given in the case files referenced in section \ref{sec:data}.

\begin{table}[!tbp]
\renewcommand{\arraystretch}{1.3}
\caption{Selection of relevant simulatin parameters of the Laser Powder Bed Fusion test case. The subscripts l, s and g signify the liquid, solid and gaseous phases of 316L stainless steel involved in the model. \cite{valenciaThermophysicalProperties2008,brooksSurfaceTensionSteels2005}}
\label{tab:lpbf-parameters}
\centering
\begin{tabular}{ll}
\toprule
Symbol & Value [Unit]\\
\midrule

  Laser Power & \SI{300}{\watt} \\
  Laser scan speed & \SI{0.9}{\metre \per \second} \\
  Density $\rho_{s}$ & \SI{7950}{\kilo\gram \per \metre\cubed} \\
  Heat capacity $c_{p,s}$ & \SI{412}{\joule \per \kilo\gram \per \kelvin} \\
  Thermal conductivity $\lambda_{s}$ & \SI{1.2e-6}{\watt \per \metre \per \kelvin} \\
  Enthalpy of fusion $H_{f}$ & \SI{2.6e5}{\joule \per \kilo\gram} \\
  Density $\rho_{l}$ & \SI{6881}{\kilo\gram \per \metre\cubed} \\
  Heat capacity $c_{p,l}$ & \SI{790}{\joule \per \kilo\gram \per \kelvin} \\
  Viscosity $\mu_{l}$ & \SI{5.85e-3}{\kilo\gram \per \metre \per \second\squared} \\
  Thermal conductivity $\lambda_{l}$ & \SI{6.6}{\watt \per \metre \per \kelvin} \\
  Heat capacity $c_{p,g}$ & \SI{1900}{\joule \per \kilo\gram \per \kelvin} \\
  Viscosity $\mu_{g}$ & \SI{1.8e-5}{\kilo\gram \per \metre \per \second\squared} \\
  Thermal conductivity $\lambda_{g}$ & \SI{6.6}{\watt \per \metre \per \kelvin} \\
Surface tension coefficient $\sigma$ & \SI{1.908}{\newton \per \metre} \\
Marangoni coefficient $\sigma_{T}$ & \SI{-1.622e-5}{\newton \per \metre \per \kelvin} \\

\bottomrule
\end{tabular}
\end{table}

\begin{figure}[!t]
  \centering
  \includegraphics[width=\linewidth]{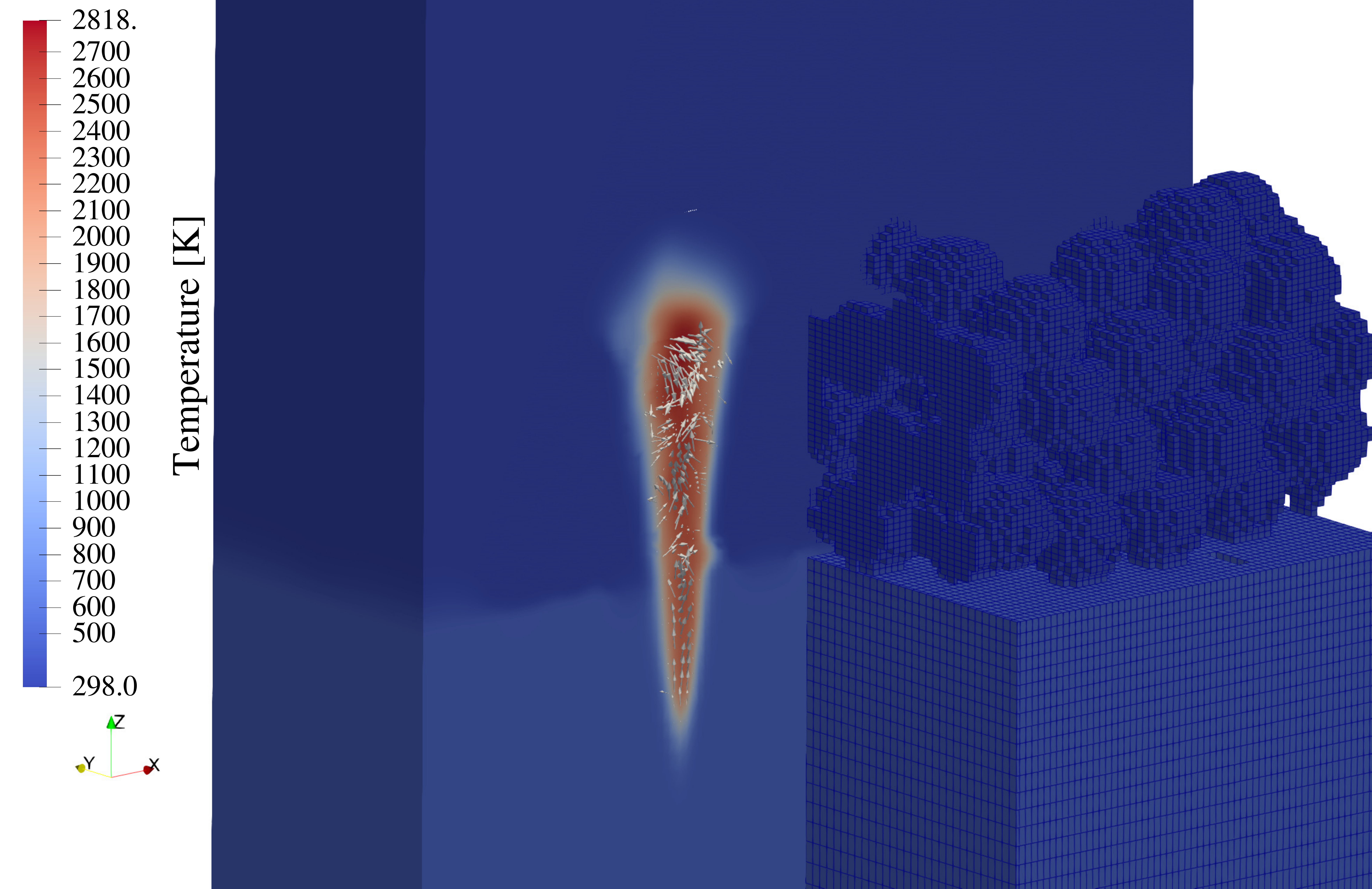}
  \caption{Temperature field of the investigated Laser Powder Bed Fusion process. The velocity field of the melt pool is indicated by arrows. Additionally, the discretised powder bed is shown along with the computational grid as a reference.}
  \label{fig:lpbf-benchmark}
\end{figure}

In Fig. \ref{fig:lpbf-benchmark}, we can observe the complex flow field that evolves within the molten material. The melt region largely penetrates the already solidified material below the powder bed, which agrees well with experimental observations done on laser melting of 316L stainless steel \cite{mohrInSituDefectDetection2020}.
Examining the cross-section further, we observe high flow velocities at the surface of the molten material, indicated by arrows in Fig. \ref{fig:lpbf-benchmark}. The large spatial temperature gradients produce variations in surface tension and hence lead to Marangoni convection at the surface. This current of liquid metal helps to distribute the heat that is input by the laser over the entire melt pool.
Furthermore, we can examine a small eddy forming within the melt pool, indicated by the circularly arranged arrows in Fig. \ref{fig:lpbf-benchmark}. This is also to be expected, albeit not due to Marangoni convection, and assists the homogenization of the temperature field within the melt pool.

\section{Summary and Future Work}

In this work, we proposed a numerical solver based on the Finite Volume method that is able to incorporate various physics that are relevant for modelling non-isothermal surface flows of arbitrary many phases. We demonstrated the capabilities of the framework using one simple, two-dimensional case showing the effects of marangoni convection and subsequently with a considerably more complex and application-driven process modelling a scenario common in additive manufacturing.

We expect that the results of this work will help to improve understanding of the complex physics of non-isothermal surface tension driven flows. Furthermore, the high fidelity data that can be generated using the proposed solver can be used in order to train physics informed machine learning models \cite{karniadakisPhysicsinformedMachineLearning2021}. It has previously been shown that such models can speed up simulations by order of magnitudes by execution of a simple forward pass. However, accurate and rich training data are needed in order to train such models \cite{liFourierNeuralOperator2021,luDeepONetLearningNonlinear2020}. We anticipate that this work can help supply this data and further enhance the field of modelling for complex manufacturing processes, as experimental data for the applications and regimes discussed are not easy to generate.

Further means of extending the proposed framework include the incorporation of ray tracing algorithms in order to accurately predict heat source dynamics, further tweaking of the underlying volume of fluid formulation for robustness and accuracy and improving parallel solving capabilities on specialized hardware, such as General Purpose GPU architectures (GPGPU).

\section*{Acknowledgment}

We thank Martin Kronbichler and Peter Munch from University of Augsburg, Germany for valuable discussions and input on the presented topic.\\
Magdalena Schreter is supported by the Austrian Science Fund (FWF) via a FWF Schrödinger scholarship, FWF project number J-4577-N.

\section*{Data Availability}\label{sec:data}

The presented solver is available via Github at \href{https://github.com/pzimbrod/thermocapillaryInterFoam}{pzimbrod/thermocapillaryInterFoam}. All data used in this article along with the case files discussed in sections \ref{sec:droplet} and \ref{sec:lpbf} can be accessed at \href{https://github.com/pzimbrod/ieee-iceccme-2022}{pzimbrod/ieee-iceccme-2022}.

\IEEEtriggeratref{23}
\bibliographystyle{ieeetr}
\bibliography{references}

\end{document}